\documentstyle[amsfonts,12pt]{article}

\def\half{\textstyle{\frac{1}{2}}}

\def\at{{\tilde a}}
\def\vt{{\tilde v}}

\def\ra{\rightarrow}
\def\tint{{\textstyle\int}}

\def\s{\hskip.08em}

\def\om{\omega}
\def\a{\alpha}
\def\b{\begin{eqnarray*}}     
\def\e{\end{eqnarray*}}       
\def\bn{\begin{eqnarray}}     
\def\en{\end{eqnarray}}       
\def\<{\langle}
\def\>{\rangle}

\def\no{\nonumber}

\def\{{\lbrace}
\def\}{\rbrace}
\begin{document}
\title{Signal Transmission in Passive Media}
\author{John R. Klauder}
\date{}     
\maketitle
\begin{abstract}
{\bf Abstract:} Under rather general assumptions, and in a relatively simple and straightforward manner, it is shown that the characteristics of signals which travel through homogeneous, as well as inhomogeneous,  passive media have the principal features usually associated with the phenomena of precursors, as generally follows from more elaborate studies. The simplicity of the present arguments permit analytic studies to be made for a greater variety of media than is normally the case. 
\end{abstract}  

\section*{{\bf 1. Introduction}}
The propagation of electromagnetic signals through passive media, i.e., media exhibiting both dispersive and absorptive properties,  of many different kinds is a problem of interest in many applications. For simplicity of analysis, we model such problems in relatively economical but fairly general terms. Let $f_0(t)$ denote the temporal dependence of the initial signal at some initial depth level (say, $z=0$) inside, or at the surface, of some medium. We presume the medium to be homogeneous at this stage and that its response characteristics are stationary. We denote by $f_z(t)$ the temporal dependence of the signal after traveling a distance $z>0$ through the passive medium. The relation between the input and output signals is that of a standard filtration which may be presented in several different forms. One such form is given by convolution with the medium impulse response function $m_z(t)$, and takes the form           
  \bn f_z(t)=\int m_z(t-s)\s f_0(s)\,ds\;.  \en
Note that all integrals with unspecified limits are assumed to extend over the whole real line. As a convolution, these functions have a multiplicative relation in frequency space given by
  \bn F_z(\om) =M_z(\om)\s F_0(\om)\;, \en
where 
  \bn && \hskip.1cm F_z(\om)=\int e^{i\om t}\,f_z(t)\,dt\;, \\
      && M_z(\om)=\int e^{i\om t}\,m_z(t)\,dt\;,  \en
denote the corresponding Fourier transforms.
As a causal phenomena (no signal out before a signal in), it follows that $m_z(t)\equiv0$, for $t<0$, a property that means that $M_z(\om)$ has no poles in the upper half complex $\om$ plane, i.e., no poles for $\om =\om_r+i\s\om_i$ with $\om_i>0$.  These few properties are standard for causal signal transmission \cite{jac}.

Let us now take up an elementary property of $M_z(\om)$ regarding its dependence on the depth parameter $z$. Clearly,
  \bn &&F_{z_1+z_2}(\om)= M_{z_1}(\om)\s F_{z_2}(\om)\no\\
         &&\hskip1.8cm =M_{z_1}(\om)\s M_{z_2}(\om)\s F_0(\om)\no\\
   &&\hskip1.8cm =M_{z_1+z_2}(\om)\s F_0(\om)\;. \en
Since $F_0(\om)$ is arbitrary in this expression, it follows that
  \bn M_{z_1+z_2}(\om)=M_{z_1}(\om)\s M_{z_2}(\om)  \label{e7}\en
holds for all nonnegative $z_1$ and $z_2$. In particular, it follows that 
  \bn M_0(\om)=M_0(\om)^2  \en
implying that either $M_0(\om)=1$ or $M_0(\om)=0$. If there are frequency bands that $F_0(\om)$ will never possess, then it is possible to choose $M_0(\om)=0$ in that band, but for simplicity we may as well assume that $M_0(\om)=1$ for all $\om$ since the final answer will be unchanged. Moreover, on physical grounds it is natural to assume, for each $\om$, that $M_z(\om)$ is a continuous function of $z$. In that case, a standard argument asserts that the only function that fulfills (\ref{e7}) necessarily has the functional form 
  \bn  M_z(\om)= e^{-z\s A(\om)} \;. \label{e33}\en
Let us recall that argument.

To establish the stated form for $M_z(\om)$, we first observe that
  \bn M_z(\om)=[\s M_{z/n}(\om)\s]^n   \en
for an arbitrary positive integer $n$. Since $M_z(\om)$ is assumed to be continuous and $M_0(\om)=1$, it follows that $M_{z/n}(\om)\neq0$ for large enough
$n$, and so $M_z(\om)\neq0$ for any $z$. Next, let $R_z(\om)\equiv \ln(M_z(\om))$, which is well defined for all $z$ and $\om$. Furthermore, (\ref{e7})
leads to $R_{z_1+z_2}(\om)=R_{z_1}(\om)+R_{z_2}(\om)$. This relation implies that $R_{nz}(\om)=n\s R_z(\om)$, as well as $R_{z/m}(\om)=(1/m)\s R_z(\om)$,
for arbitrary positive integers $n$ and $m$. Consequently,
   \bn  R_{nz/m}(\om)=(n/m)\s R_z(\om)\;.   \en
Now let $z=1$ so that $R_{n/m}(\om)=(n/m)\s R_1(\om)$; then take a limit such that $(n/m)\ra z$, an arbitrary positive number. By the assumption of continuity, this limit converges and we learn that
  \bn R_z(\om)=z\s R_1(\om)\equiv -z\s A(\om)\;.  \en
This concludes the proof of (\ref{e33}).

As the transfer function of a passive medium it is necessary that $|M_z(\om)|\le1$, in which case we also find that ${\rm Re}\s A(\om)\ge0$. To simplify notation, we let
  \bn A(\om)=B(\om)+i\s C(\om)\;, \en
where $B(\om)={\rm Re}\s A(\om)$ and $C(\om)={\rm Im}\s A(\om)$. 
Thus, we require that  
   \bn B(\om)\ge0 \;. \en
By causality, $B(\om)$ and $C(\om)$ are related by Hilbert transforms, while $C(\om)$ is undetermined by this procedure up to the contribution of a so-called Blaschke factor \cite{blas}.
 We also note that, for real functions $f_z(t)$ for all $z$ and any initial signal $f_0(t)$, it necessarily follows that the impulse response function $m_z(t)$ is a real function. Consequently, $M_z(-\om)=M_z(\om)^*$, which implies that
  \bn  && B(-\om)=B(\om) \;, \\ 
       && C(-\om)=-C(\om)\;.  \en
Since $ B(\om)\ge 0$, it follows that the function $B(\om)$ has special properties at and near any frequency, say $\om'$, for which $B(\om')=0$. Assuming on physical grounds that $B(\om)$ is continuous in $\om$, and even has a few derivatives [e.g., $B(\om)\in C^2$], then the behavior of $B$ in the vicinity of $\om'$ is given by
  \bn B(\om)=\half a^{-1}\s(\om-\om')^2\;, \hskip1cm (\om-\om')^2\ll1\;. \en
where $a>0$ has the units of acceleration (e.g., $m/s^2$). 

Let us focus on the behavior of $A(\om)$ in the vicinity of $\om=0$. In that case, and to sufficient accuracy, we have 
  \bn && B(\om)=\ell^{-1}+\half\s a^{-1}\s \om^2\;, \label{f333}\\
      && C(\om)=-v^{-1} \om\;,  \label{f334}\en
where $\ell$, $v$, and $a$ denote constants with the dimensions of length, velocity, and acceleration, respectively. The behavior at $\om=0$ forces $\ell^{-1}\ge0$. Moreover, when $B(0)=\ell^{-1}=0$, it follows that $a^{-1}\ge0$; in that case, we shall adopt the generic form and assume that $a>0$. 

It is natural to ask the question: Under what circumstances does it follow that $B(0)=0$? A common model of a passive medium is based on the Lorentz model \cite{lor}. The physical principle behind this model is the response of a damped harmonic oscillator to a driving force. The physical model in this case is given by
  \bn m\s{\ddot q}(t)+b\s{\dot q}(t)+k\s q(t)=h(t)\;, \en
where $m>0$, $b\ge0$, and $k>0$ denote the inertial, damping, and restorative parameters, $q(t)$ denotes the oscillator amplitude, and $h(t)$ denotes the applied force. Assuming that any transients have decayed away, let us focus on the steady state solution. If $h(t)=U\cos(\om\s t)$, then $q(t)=V\cos(\om\s t)+W\sin(\om\s t)$, where
  \bn && \hskip.06cm V=\frac{(k-m\s\om^2)\s U}{(k-m\s\om^2)^2+(b\s\om)^2}\;,\\
  &&   W=\frac{b\s\om\s U}{(k-m\s\om^2)^2+(b\s\om)^2}\;.  \en
If the driving frequency $\om=0$ for this model it follows that $V=U/k$ and $W=0$. Since this result is entirely independent of the dissipation coefficient $b$, it holds even when $b=0$, namely, it holds in the case there is no loss whatsoever. Consequently, for such a model it follows that at $\om=0$ there can be {\it no} loss and so $B(0)=0$. Of course, such an argument does not apply
to a metal, but it does hold for a broad spectrum of materials.

Although we have argued that $B(0)=0$ for a specialized model, it is quite plausible that such a result holds for a wide class of models as well. One such more general model would be of the form
  \bn m_i\s{\ddot q}_i(t)+b_i\s{\dot q}_i(t)+\Sigma_{j=1}^J\s k_{ij}\, q_j(t)=h(t)\;, \hskip1cm 1\le i\le J\;, \en
which corresponds to a situation where $h(t)$ acts as the driving force on a collection of coupled damped harmonic oscillators. Again we may consider the steady state response when $h(t)=U\cos(\om\s t)$. In particular, if $\om=0$ the solution is totally independent of any dissipation parameters. Consequently, it follows that $B(0)=0$ in this case as well. It is not difficult to imagine other, more complicated situations for which the same result holds. Therefore, on the basis of a wide class of physically plausible models, we are encouraged to consider the consequences of the assumption that $B(\om)=\half\s a^{-1}\s\om^2$, $a>0$, when $\om$ is small.

We now examine the signal behavior as a function of $z$ that follows from the foregoing discussion. In particular, we first focus on the impulse response function
  \bn m_z(t)=\frac{1}{2\pi}\int e^{-i\om t}\,e^{iv^{-1}z\s\om- za^{-1}\om^2/2-z\s {\tilde A}(\om)}\; d\om\;,  \en
where ${\tilde A}(\om)\equiv{\tilde B}(\om)+i\s{\tilde C}(\om)$, and ${\tilde B}(\om)=O(\om^4)$ while ${\tilde C}(\om)=O(\om^3)$. To simplify matters, let us further assume that ${\tilde B}(\om)$ is bounded away from $0$ for nonzero values of $\om$. In that case, for large $z$, we observe that the integrand tends to be strongly suppressed for $\om$ away from zero. To sufficient accuracy it follows that we may set
  \bn && m_z(t)=\frac{1}{2\pi}\int e^{-i\om(t-v^{-1} z)}\,e^{- za^{-1}\s\om^2/2}
\;d\om\no\\
 &&\hskip1.1cm =\sqrt{\frac{a}{2\pi z}}\,\exp{\!\bigg(-\frac{a(t-v^{-1}z )^2}{2z}\bigg)}\;. \label{e333}\en
Consequently, for large $z$ values we learn that
  \bn && f_z(t)=\int m_z(t-s)\s f_0(s)\,ds \no\\
  &&\hskip.95cm=\sqrt{\frac{a}{2\pi z}}\int \exp\!\bigg(-\frac{a(t-s-v^{-1} z)^2}{2z}\bigg)\, f_0(s)\,ds\no\\
  &&\hskip.95cm=\sqrt{\frac{a}{2\pi z}}\int e^{-a\s s^2/2z}\,f_0(s+t-v^{-1} z)\,ds\;. \label{e22}\en

There are two characteristic features of this resultant behavior that are worth emphasizing:\vskip.2cm
\hskip.3cm 1. According to the square root prefactor, it follows that the amplitude decays proportional to $z^{-1/2}$. \vskip.2cm
\hskip.3cm 2. According to the integral factor, the result has a temporal dependence determined by the form of the original signal centered at a
temporal value of $t-v^{-1} z$, and then integrated over that region weighted by a Gaussian of increasing width.

\subsubsection*{{\it 1.1 Approximate causality}}
Before passing to several examples of (\ref{e22}) it is important to observe that the medium impulse response function (\ref{e333}) is {\it not}
causal since it does not vanish for $t<0$. In point of fact, the expressions
(\ref{f333}) and (\ref{f334}) do {\it not} make a strictly causal impulse response function. However, they make an {\it approximate} causal function provided
that $az\gg v^2$ since in that case the impulse response function is extremely
small for $t<0$. Let us examine this situation more closely.

In order for the medium function $A(\om)$ to correspond to a causal
function it should be represented in the form
  \bn  A(\om)=\int_0^\infty[1-e^{i\om\s s}]\s k(s)\,ds\;,  \en
where $A(0)=0$, and in order that ${\rm Re}A(\om)=B(\om)\ge0$, it is sufficient that
the weight function $k(s)\ge0$. Observe, by this representation, that
$A(\om)$ is analytic in the upper half plane, which is a necessary and
sufficient condition for the impulse response function to be causal. Now, 
as an example, consider $k(s)=K' e^{-Ks}$, where $K'$ and $K$ are positive
parameters to be determined. In that case,
  \bn && B(\om)=K'\s\int_0^\infty[1-\cos(\om\s s)]\,e^{-Ks}\,ds \;,\\
  && C(\om)=-\s K'\s\int_0^\infty \sin(\om\s s)\,e^{-Ks}\,ds  \;.  \en
We are interested in a large value for both $K$ and $K'$ such
that, to some approximation, $k(s)$ acts rather like a delta
function in the integrand focusing on values of $s$ near zero. Thus, to a
leading approximation,
  \bn  &&B(\om)={\half}\s K'\s\om^2 \int_0^\infty s^2 \s 
e^{-Ks}\,ds +\cdots \;, \label{g30}\\
  &&C(\om)= -\s K'\s\om\int_0^\infty s \,e^{-Ks}\,ds+\cdots\;. \label{g31} \en
On comparison with (\ref{f333}) and (\ref{f334}) we see that
  \bn && a^{-1}=K'\int_0^\infty s^2\s e^{-Ks}\,ds = 2K'/K^3\;,  \\
    && v^{-1}=K'\s\int_0^\infty s\,e^{-Ks}\, ds =K'/K^2\;,  \en
and so $a=\half\s K\s v$. Provided $K\gg1$, this approximation is
satisfactory from a practical viewpoint; otherwise, some of the additional 
terms in (\ref{g30}) and (\ref{g31}) need to be taken into account. 
Unfortunately, introducing additional
terms makes the medium impulse response function no longer analytically
tractable, generally speaking. 
 
\subsubsection*{{\it 1.2 Specific examples}}
Let us illustrate the resultant signal (\ref{e22}) for two important examples. 
First we consider an analytically simple case, namely a Gaussian initial pulse. Thus, 
let
  \bn  f_0(t)=e^{-t^2/2T^2}\s\cos(\om_o t) \;, \en
where we refer to $T$ as the pulse width.
In this case, (\ref{e22}) reads
  \bn f_z(t)=\sqrt{\frac{a}{2\pi z}}\int e^{-a\s (s-t+v^{-1}z)^2/2\s z}\,e^{-s^2/2\s T^2}\,\cos(\om_o s)\;ds\;. \en
An elementary integration leads to 
 \bn &&f_z(t)=\frac{1}{\sqrt{1+z/a\s T^2\s}}\no\\
&&\hskip1cm\times\exp\!\bigg(-\frac{[(a/z)(t-v^{-1}z)^2+\om_o\s T^2]}{2(1+a\s T^2/z)}\bigg)\s\cos\bigg(\frac{\om_o(t-v^{-1}z)}{1+z/a\s T^2}\bigg)\;. \label{e35}\en
If we assume that $z\gg a\s T^2$, then in effect (\ref{e35}) reduces
to 
  \bn f_z(t)=\sqrt{\frac{a\s T^2}{z}}\,\exp\!\bigg[-\frac{a\s(t-v^{-1}z)^2}{2\s z}-\frac{(\om_o T)^2}{2}\bigg]\;. \en
The peak amplitude of the resultant pulse $f_z(t)$ occurs when $z=v\s t$, or, stated otherwise when $t=v^{-1}z$. At that point in $z$ (or at that moment in time $t$), the amplitude has been reduced from its initial maximum of unity by a factor proportional to $z^{-1/2}$. The resultant pulse shape in the present case is still Gaussian, but it is substantially broadened since the pulse width $T$ of the input Gaussian has been replaced by $\sqrt{z/a}\gg T$. If we assume the energy (or power) in a given pulse $f_z(t)$ is proportional to $\tint|f_z(t)|^2\,dt$, then it follows, for $z\gg a\s T^2$,  that the normalized signal energy is given by
  \bn  \frac{\tint|f_z(t)|^2\,dt}{\tint|f_0(t)|^2\,dt}=2\sqrt{\frac{a\s T^2}{z}}\,\frac{e^{-(\om_o T)^2}}{1+e^{-(\om_o T)^2}}\;. \en
As anticipated, this ratio decreases as $z$ increases. However, thanks to the pulse broadening, it decreases only as $z^{-1/2}$ rather than $z^{-1}$, as may have been expected on the basis of (\ref{e22}).  

As a second example, we suppose that
  \bn f_0(t)={\rm rect}(t/T)\s\cos(\om_o \s t) \;, \label{i2}\en
which leads to the output signal
  \bn f_z(t)=\sqrt{\frac{a}{2\pi z}}\int_{-T/2}^{T/2} e^{-a(s-t+v^{-1}z)^2/2z} \;\cos(\om_o\s s)\,ds\;. \label{e44}\en
This integral depends on $z$ in a moderately involved way, and for convenience 
we only evaluate (\ref{e44}) in the case that $z\gg aT^2$. Then we find
  \bn &&f_z(t)=\sqrt{\frac{a}{2\pi z}}\int_{-T/2}^{T/2} e^{-a(t-v^{-1}z)^2/2z} \;\cos(\om_o\s s)\,ds\no\\
  &&\hskip1cm=\sqrt{\frac{a\s T^2}{2\pi z}}\,\bigg[\frac{\sin(\om_o\s T/2)}{(\om_o T/2)}\bigg]\,e^{-a(t-v^{-1}z)^2/2z}\;. \label{g40}\en
Again, the
result is that of a traveling pulse moving essentially at the speed $v$, the
amplitude of which falls off with distance like $z^{-1/2}$, provided of course that $\sin(\om_o T/2)\neq0$. For the case that $\sin(\om_o T/2)=0$, $\om_o T\neq0$, see below.

\subsubsection*{{\it 1.3 More general inputs}}
Apart from an unimportant coefficient, observe that the amplitude factor arising from $f_0(t)$ in (36), i.e., $\exp(-\om_o^2T^2/2)$, is basically $F_0(0)=\int f_0(s)\, ds$, and the amplitude factor in (\ref{g40}), i.e., $\sin(\om_o T/2)/(\om_oT/2)$ is essentially just $F_0(0)=\int f_0(s)\s ds$. We can understand and generalize this result as follows.  Since
  \bn f_z(t)=\tint m_z(t-s)\s f_0(s)\,ds\;,  \en
it follows, for an initial signal which is much narrower in the temporal domain
than the medium response function, that we can make a Taylor series expansion
of the medium response function. In particular,
we have  
  \bn  &&f_z(t)=\tint [\s m_z(t)\s f_0(s)-m_z'(t)\s sf_0(s)+
\half m_z''(t)\s s^2f_0(s) + \cdots\,]\,ds \no\\
  &&\hskip.98cm =m_z(t)\tint f_0(s)\s ds-m_z'(t)\tint sf_0(s)\s ds
+\half m_z''(t)\s \tint s^2f_0(s)\s ds +\cdots \;. \no\\
&&  \label{i3}\en
When $m_z(t)$ is sufficiently broad in time, then its time derivatives $m'_z(t),\,m''_z(t)$, etc., are 
typically much smaller. In this case, the indicated sum involves terms that are
rapidly diminishing relative to the former terms.

As an example of the application of (\ref{i3}), suppose that $\om_oT=2n\pi$,
$n\in\{1,2,3,\ldots\}$, then (\ref{g40}) vanishes because $\tint f_0(s)\s ds=0$. Moreover, since $f_0(t)$ in (\ref{i2}) is an even function, it also follows
that $\tint s\s f_0(s)\s ds=0$. Thus, according to (\ref{i3}), we find to
leading order that
  \bn f_z(t)=\half\s m''_z(t)\,\tint s^2\s f_0(s)\,ds \;,  \en
namely,
  \bn f_z(t)=\frac{(-1)^n}{2(n\pi)^2}\s\sqrt{\frac{a\s T^2}{2\pi z}}\s\bigg(\frac{a\s T^2}{z^2}\bigg)\s\bigg[\s a(t-z/v)^2-z\s\bigg]\,e^{-a(t-z/v)^2/2z} \;,  \en
where the amplitude has been evaluated when $\om_oT=2n\pi$.

\subsubsection*{{\it 1.4 Commentary}}
Observe, that although the original Gaussian and the rectangular envelopes are fundamentally different, the pulses that result after a long propagation distance $z\gg a\s T^2$ are essentially the same, differing only in their respective amplitudes. Those amplitudes effectively measure the spectral content of the original pulse at zero frequency where the narrow window of transmission exists, and for the Gaussian envelope that content is much smaller than is available with the rectangular envelope. This result is entirely consistent with the concept that the output selects its preferred signal shape consistent with the frequency window allowed by the medium. Ultimately, this means that the transmitted signal is Gaussian in character with an ever increasing pulse width in accord with the ever narrowing Gaussian pass band as $z$ increases. 
This kind of temporal behavior is essentially the kind of behavior typically ascribed to the so-called Brillouin precursor \cite{lor}. It is paradoxical, but true, that Brillouin precursors require pulse 
envelopes with {\it high} frequency content. Sharp edges in pulse envelopes 
significantly enhance the high frequency content as compared to smooth 
envelopes, and this situation is clearly illustrated in Eqs.~(36) and (40). 

However, sharp edges in pulse envelopes are not the only way to raise the 
high frequency content.

\section*{{\bf 2. Chirp signals}}
As an initial wave form consider the linearly chirped signal \cite{chi}
  \bn  f_0(t)=e_T(t)\s \cos(\om_o\s t+\half\a t^2)    \label{t3} \;,  \en
where the signal envelope $e_T(t)$ has a pulse width of $T$ and we assume 
that 
$\a T^2\gg1$. The Fourier transform of (\ref{t3}) may be approximately 
given by   
\bn &&F_0(\om)=\frac{1}{2}\s\sqrt{\frac{2\pi}{i\s\a}}\s 
e_T(\s(\om -\om_o)/\a)\s
e^{i(\om-\om_o)^2/2\a}\no\\
&&\hskip1.4cm +\frac{1}{2}\s\sqrt{\frac{2\pi\s i}{\a}}\s e_T(-(\om +\om_o)/\a)
\s e^{-i(\om+\om_o)^2/2\a}\;.  \en
The frequency content at $\om=0$ is then given by
  \bn F_0(0)=\frac{1}{2}\s\sqrt{\frac{2\pi}{i\s\a}}\s e_T(-\om_o/\a)\s
e^{i\s\om_o^2/2\a}+\frac{1}{2}\s\sqrt{\frac{2\pi\s i}{\a}}\s e_T(-\om_o/\a)
\s e^{-i\s\om_o^2/2\a}\;.  \en
As an illustration, consider $e_T(t)=\exp(-t^2/2T^2)$. In this case,
 \bn F_0(0)=\sqrt{\frac{\pi}{\a}}\s \exp\bigg[-\frac{(\om_oT)^2}{2\a^2 T^4}
\bigg]\s\bigg[\s\cos\bigg(\frac{(\om_oT)^2}{\a T^2}\bigg)+\sin\bigg(
\frac{(\om_oT)^2}{\a T^2}\bigg)\bigg]\;.  \en
As noted above, this factor essentially dictates the amplitude for the output Gaussian in (24).

The main point to observe in this calculation is the significant 
enhancement of the
overall amplitude of the signal that will emerge. In particular, the 
primary change has been that the factor $\exp[-(\om_oT)^2/2]$ in (36) has been 
replaced by $\exp[-(\om_oT)^2/2\a^2T^4]$, where $\a T^2\gg1$. This model calculation 
strongly suggests that combining chirp signals with envelopes having sharp 
or nearly sharp edges should lead to
signals with relatively significant precursor content. 

\section*{{\bf 3. Composite media}}
Up to this point we have assumed that the signal propagates through a single,
homogeneous medium, the properties of which are characterized by the complex medium function $A(\om)$. Suppose instead we deal with a composite media made
up from $J$ layers, each of thickness $l_1=z_1,\;l_2=z_2-z_1,\ldots,l_J=
z_J-z_{J-1}$. To each layer we attribute a medium function $A_j(\om)$, 
$1\le j\le J$. If we assume that $z_j\le z\le z_{j+1}$, for example, then
  \bn  M_z(\om)=\exp[-(z-z_j)\s A_{j+1}(\om)-\Sigma_{k=1}^j\,l_kA_k(\om)]\;. \en
A more convenient expression for such media is given by
  \bn  M_z(\om)=\exp[-\tint_0^zA_u(\om)\,du]\;,  \en
which applies for piecewise-constant functions $A_u(\om)$. It is plausible
that such an expression applies as well for slowly varying media functions
$A_u(\om)$. 

As was previously the case, it is useful to separate $A_u(\om)$ into real
and imaginary parts, i.e.,
  \bn  M_z(\om)=\exp[-\tint_0^zB_u(\om)\,du-i\tint_0^zC_u(\om)\,du]\;.  
\label{f8}\en
Under suitable conditions we may set
\bn  B_u(\om)=\half a^{-1}(u)\s\om^2\;,\hskip1cm C_u(\om)=-v^{-1}(u)\s\om\;,\en
in which case we are led to
  \bn M_z(\om)=\exp[-\half\tint_0^za^{-1}(u)\s\om^2\,du+
i\tint_0^zv^{-1}(u)\s\om\,du\s]\;. \label{f9}\en
Equations (\ref{f8}) and (\ref{f9}) apply to various situations.

As a simple example of how these equations may be applied let us consider a signal that penetrates a passive medium with a finite thickness of $\ell$ and
otherwise propagates in free space. Therefore, in this case $B_u(\om)>0$
for $0<u<\ell$, while for $u>\ell$, $B_u(\om)=0$. This
example pertains to a radar signal that penetrates a canopy of foliage or some other absorbing and dispersive medium. Even a comparatively thin, passive medium
may be well approximated by (\ref{f9}) when $A_u(\om)$ is sufficiently large for $\om\neq0$. Under the stated conditions it follows, when $z>\ell$, that
  \bn  M_z(\om)=\exp\{-\ell\om^2/2\at +i[(z-\ell)\om/c+\ell\om/\vt]\} \;, \en
where $c$ is the speed of light in vacuum, and
  \bn && \at\equiv \ell\,[\tint_0^\ell \s a(u)^{-1}\,du]^{-1}\;, \no\\
  &&\vt\equiv\ell\,[\tint_0^\ell \s v(u)^{-1}\,du]^{-1}\;.  \en
Since $v(u)\le c$, it follows that $\vt\le c$. For simplicity in our example, let us assume that $\vt\simeq c$. Finally, when $\ell/\at$ is
sufficiently large, we are led to 
  \bn  f_z(t)=\sqrt{{\at}/{\ell}}\; e^{-\at(t-z/c)^2/2\ell}\,
\tint f_0(s)\,ds\;, \en
defined for $z>\ell$. If $F_0(0)=0$, then appeal can be made to (42).

This equation should be compared with Eq.~(36) that corresponds to propagation within a single, homogeneous passive medium. The temporal signal in both cases is Gaussian and apart from the overall amplitude, the shape is principally determined by the passive medium. This result holds whenever the medium acts 
as a low
pass filter for an initial signal that, comparatively speaking, has a broad 
spectral content including $\om=0$.

\section*{{\bf 4. Stochastic media}}
Under certain circumstances, it is appropriate to regard the passive medium as a stochastic variable and to define the output signal as the ensemble average of a stochastic output. A description of this kind may arise if the output signal is the superposition of multiple signals from a rapidly changing collection of media profiles. More specifically, such a situation may arise in the case
of a coherent superposition of multiple radar pulse returns that have traversed a rapidly fluctuating medium. While details of the fluctuating medium may be too difficult to describe, certain overall features -- stochastic invariants -- may provide gross characteristics with which selected signal properties may be analyzed.

Let us focus on propagation through a stochastic but otherwise homogeneous medium. In mathematical terms this means that $M_z(\om)$, and thereby
$A(\om)$, are regarded as suitable stochastic variables. The observed signal $f^o_z(t)$ is determined by the expression
  \bn  f^o_z(t)=\frac{1}{2\pi}\int e^{-i\om t}\,\<\s e^{-zA(\om)}\s\>\s F_0(\om)
\,d\om\;,  \en
where $\<(\,\cdot\,)\>$ denotes an ensemble average over medium characteristics. Let us illustrate the kind of behavior that may occur. For that purpose we again suppose that 
  \bn A(\om)=\half a^{-1}\om^2-iv^{-1}\om  \en
and we assume that $a$ is the only random variable. As an example we suppose that the variable $a$, which can only be positive, is distributed according to a probability density given by
   \bn p(a^{-1})=\frac{b}{m!}\s(ba^{-1})^m\s e^{-ba^{-1}}\;,
\hskip1cm 0<a^{-1}<\infty\;,  \en
where $b$ is a parameter which essentially sets the scale. In particular,
the mean value of $a^{-1}$ is given by
  \bn \<\s a^{-1}\s\>=\int_0^\infty a^{-1}\s p(a^{-1})\, da^{-1} =\frac{m+1}{b} \;. \en
As a consequence, we find that
   \bn  f^o_z(t)=\frac{1}{2\pi}\int e^{-i\om(t-z/v)}\frac{F_0(\om)}
{(1+z\s\om^2/b)^{m+1}}\,d\om  \;.  \label{i77}\en

Even before evaluating the Fourier transform it is clear that our stochastic example has led to a dramatically enhanced pass band at low frequencies. This result is due, of course, to the nonvanishing probability that arbitrarily 
small values of $a^{-1}$ are present in the ensemble, despite the fact that
 the probability that $a^{-1}=0$ is strictly zero for any $m\ge1$. 

As outlined in (\ref{i3}), the leading behavior of (\ref{i77}) is given by
  \bn  f^o_z(t)=m^o_z(t)\s \tint f_0(s)\,ds\;,  \en
where
   \bn  m^o_z(t)=\frac{1}{2\pi}\int \frac{e^{-i\om t}}{(1+z\om^2/b)^{m+1}}
\,d\om\;.  \en
This expression may be evaluated in the following form: 
  \bn  m^o_z(t)=\frac{1}{2^{m+1}}\,\frac{1}{m!}\,\sqrt{b/z}\,\bigg\{\sum_{l=0}^m\,C^{(m)}_l\,\bigg(\sqrt{b/z}\;|t-z/v|\bigg)^l\bigg\}\,e^{-\sqrt{b/z}\,|t-z/v|}\;,  \en
where the $C^{(m)}_l$ coefficients, for $0\le l\le m$ and $0\le m$, are determined by the relations
\bn && C^{(m)}_m=1 \;,\no\\
    && C^{(m)}_0=(2m-1)!!  \;,\no\\
   && C^{(m)}_l=(2m-1-l)\s C^{(m-1)}_l+C^{(m-1)}_{l-1},\hskip.5cm 1\le l\le m-1 \;. \en

\section*{{\bf 5. Conclusions}}
Traditional radar applications involve signal generation, transmission, and reception, with little or no spectral modification save for a possible spectral shift (or dilation) due to relative motion. By contrast, radar signals, or electromagnetic waves more generally, that traverse passive media with normal absorptive and dispersive properties, typically suffer significant spectral distortion. For example, over a large spectral domain, an exponential attenuation effectively removes any signal energy. An exception to exponential attenuation occurs at zero frequency, when attenuation is frequently absent. 
On the basis of generally applicable arguments, and under a wide set of circumstances, it has been clearly demonstrated in this paper that the amplitude of the transmitted signal decays as the inverse square root of the distance traveled rather than exponentially. Such behavior is exactly that characteristic of so-called Brillouin precursors \cite{jac,lor}. While such phenomena are frequently discussed in terms of more advanced mathematical techniques, it is demonstrated in the present work that they can be largely understood in terms of fairly simple pass band arguments.\footnote{It is noteworthy that my colleague Tim Olson has devised quite another explanation for precursor behavior. See his forthcoming paper ``Derivation of Precursors via Finite Toeplitz Forms''.}

The simplicity of the present discussion has enabled us to determine the form of signals that have traversed inhomogeneous or stochastic media as well. The robustness of the characteristic precursor decay rate has also been established for these more general media.

\section*{{\bf 6. Acknowledgments}}
Thanks are expressed to T. Olson and S.V. Shabanov for their interest and their comments. Additionally, W. Bomstad is acknowledged for his detailed comments on the paper. This work has been partially supported by US AFOSR Grant F49620-01-1-0473.

\end{document}